\documentclass[a4paper,11pt,twoside]{amsart}
\usepackage{amsmath,amsthm,amsfonts,amssymb,latexsym,amscd}

\def\N{\mathbb N}

\def\R{\mathbb R}
\def\Z{\mathbb Z}

\def\cC{\mathcal C}

\def\ra{\rightarrow}

\def\lra{\longrightarrow}

\newcommand{\into}{\hookrightarrow}
\newcommand{\onto}{\twoheadrightarrow}

\def\={\overset{\text{def}}=}
\def\doublearrow#1#2
{\quad\raisebox{.1cm}{$\overset{#1}\lra$}  
\hskip -.67cm
\raisebox{-.1cm}{$\underset {#2}\lra$}\quad}

\def\im{\operatorname{im}}

\def\ra{\rightarrow}

\def\pt{\operatorname{point}}

\newtheorem{prop}{Proposition}[section]
\newtheorem{proprecall}{Proposition 2.11}
\newtheorem{lem}[prop]{Lemma}
\newtheorem{rrem}[prop]{Remark}

\newtheorem{eex}[prop]{Example}
\newtheorem{ddefn}[prop]{Definition}
\newtheorem{theorem}[prop]{Theorem}
\newtheorem{cor}[prop]{Corollary}

\newtheorem{hhome}[prop]{Homework}
\newtheorem{nnumber}[prop]{}
\newenvironment{num}{\begin{nnumber}\rm}{\end{nnumber}}

\newenvironment{defn}{\begin{ddefn}\rm}{\end{ddefn}}

{\catcode`@=11\global\let\c@equation=\c@prop}

\def\psc{positive scalar curvature}
\def\pscm{positive scalar curvature metric}

\def\GLR{Gromov-Lawson-Rosenberg}

\def\lra{\longrightarrow}

\def\nb-{\nobreakdash-}

\begin{document}

\title{Remarks on a conjecture of Gromov and Lawson}
\date{July 26, 2002}
\thanks{The first and the third author were supported by NSF grants;
  the second and the
third author were supported by the Max-Planck-Institute in Bonn, Germany}

\author{William Dwyer}
\author{Thomas Schick}
\author{Stephan Stolz}
\address{Dwyer and Stolz: Dept. of Mathematics\\
University of Notre Dame\\ Notre Dame, IN 46556\\ USA}
\address{Schick:
 Fachbereich Mathematik --- Universit{\"a}t G{\"o}ttingen\\
Bunsenstr.~3\\
         37073 G{\"o}ttingen, Germany}
\maketitle

\begin{abstract}
Gromov and Lawson conjectured in \cite{GL2} that a closed spin
mani\-fold $M$ of dimension $n$ with fundamental
group $\pi$ admits a \pscm\ if and only if an associated element
in $KO_n(B\pi)$ vanishes. In this note we present counter
examples to the `if' part of this conjecture for groups $\pi$
which are torsion free and whose classifying space is a manifold
with negative curvature (in the Alexandrov sense). 
\end{abstract}

\section{Introduction}
In their influential paper \cite{GL2} Gromov and Lawson
proposed the following conjecture.

\begin{num}{\bf Gromov-Lawson Conjecture
\cite{GL2}.}\label{num:GLconj}
 Let $M$ be a
smooth, compact manifold  without boundary of dimension $n\ge
5$ with fundamental group
$\pi$. Then $M$ admits a \pscm\ if and only if $p\circ
D([M,u])=0\in KO_n(B\pi)$.
\end{num}

Here $u\colon M\to B\pi$ is the map classifying the universal
covering of $M$, and $[M,u]\in \Omega_n^{spin}(B\pi)$ is the
element in the spin bordism group represented by the pair
$(M,u)$. The maps $D$ and $p$ are natural transformations
between generalized homology theories:
\[
\Omega_n^{spin}(X)\overset{D}\ra
ko_n(X)\overset{p}\ra
KO_n(X),
\]
referred to as `spin bordism', `connective $KO$\nb-homology'
and `periodic
$KO$\nb-ho\-mology', respectively. 
As the name suggests,
$KO_n(\quad)$ is  periodic in the sense that 
$KO_n(X)\cong KO_{n+8}(X)$; moreover, 
\[
KO_n(\pt)\cong 
\begin{cases}
\Z&\text{if $n\equiv 0\mod{4}$}\\
\Z/2&\text{if $n\equiv 1,2\mod{8}$}\\
0&\text{otherwise}
\end{cases}.
\]
The homology theory $ko_n(\quad)$ is  connective; i.e., 
$ko_n(\pt)=0$ for $n<0$. Moreover, the natural transformation $p$
induces an isomorphism $ko_n(\pt)\cong KO_n(\pt)$ for $n\ge 0$
(this does {\it not} hold with $\pt$ replaced by a general space
$X$). 

We remark that the  assumption $n\ge 5$ (which wasn't present
in \cite{GL2}) should be added since Seiberg-Witten
invariants show that the conjecture is false for $n=4$
even if the group $\pi $ is trivial.
In their paper
\cite{GL2} Gromov and Lawson prove the  `only if' part of
Conjecture
\ref{num:GLconj}  for some groups
$\pi$, but later Rosenberg \cite{Ro1} noticed that the  `only if'
statement does not hold e.g.\ for finite cyclic groups, since any 
lens space
$M$ admits a metric with \psc, but
$p \circ D([M,u])$ is non-trivial. He proposed the following
variant of Conjecture \ref{num:GLconj}.

\begin{num}{\bf The Gromov-Lawson-Rosenberg Conjecture
\cite{Ro3}.}
\label{num:GLRconj} 
 Let $M$ be a smooth, compact manifold  without boundary of
dimension $n\ge 5$ with fundamental group
$\pi$. Then $M$ admits a \pscm\ if and only if
$\alpha([M,u])=0\in KO_n(C^*_r\pi)$.
\end{num}

Here $\alpha$ is the following composition
\begin{equation}\label{eq:composition}
\Omega_n^{spin}(B\pi)\overset{D}\ra
ko_n(B\pi)\overset{p}\ra
KO_n(B\pi)\overset{A}\ra KO_n(C_r^*\pi),
\end{equation}
where $A$ is the  {\it assembly map}, whose
target is the  $KO$\nb-theory of the (reduced) real
group $C^*$\nb-algebra $C_r^*\pi$ (this is a norm completion
of the real group ring $\R\pi$ and is equal to the latter for
finite $\pi$).

Rosenberg proved the `only if' part of this Conjecture in
\cite{Ro2} by  interpreting the  image of 
$[M,u]\in\Omega_n^{spin}(B\pi)$ under  $\alpha$  as a
`fancy' type of index of the Dirac operator on $M$ and showing
that this index vanishes if
$M$ admits a metric of \psc. We note that if $A$ is
injective for a group $\pi$ (which can be shown for
many groups $\pi$ and is  expected to hold for
{\it all} torsion free groups according to the
Novikov-Conjecture), both conjectures are equivalent.

All partial results concerning Conjectures \ref{num:GLconj}
and
\ref{num:GLRconj} discussed so far have to do with the `only
if' part or equivalently, with finding {\it obstructions}
(like the index obstruction
$\alpha$) against the existence of \pscm s. For a few groups
$\pi$ {\it constructions} have been found of sufficiently
many \pscm s to prove the `if' part as well and
hence to prove the Gromov-Lawson-Rosenberg Conjecture
\ref{num:GLRconj}. This includes the trivial group \cite{St1},
cyclic groups and more generally all finite groups $\pi$ with
periodic cohomology \cite{BGS}. However, a counter example
has been found for $\pi=\Z/3\times \Z^4$, $n=5$
\cite{Sch}. So far,  no counter example to the
Gromov-Lawson-Rosenberg Conjecture has been found  for {\it
finite} groups $\pi$ (which actually is the  class of groups
the conjecture was originally formulated for; cf.\ 
\cite{Ro3}, Conjecture 0.1 and \S 3). 

In this paper we address the question of whether the
Gromov-Lawson Conjecture and/or Gromov-Lawson-Rosenberg
Conjecture might be true for {\it torsion free groups}. Alas,
the answer is `no' in both cases; more precisely:

\begin{theorem}\label{thm:mainthm} There are finitely
generated torsion free groups $\pi$ for which the `if'
part of the Gromov-Lawson Conjecture \ref{num:GLconj} (and
hence also the `if' part of the Gromov-Lawson-Rosenberg
Conjecture \ref{num:GLRconj}) is false. In other words,
there is a closed spin manifold $M$ of dimension $n\ge 5$ with
fundamental group $\pi$ and $p\circ D([M,u])=0$ which does not
admit a \pscm.
\end{theorem}

Still, one might hope to save these conjectures at least for
groups $\pi$ satisfying some additional geometric conditions,
say of the kind that guarantee that the assembly map $A$ is 
an isomorphism (we remark that according to the Baum-Connes
Conjecture the assembly map is an isomorphism
for all torsion free groups $\pi$). For these groups in
particular the `only if' part of the conjectures is true.

\begin{num}{\bf Examples of groups $\pi$ for which the
assembly map is an isomorphism}
\begin{enumerate}
\item Any countable group $\pi$ in the class
$\cC=\cup_{n\in
\N}\cC_n$, where the class $\cC_0$ consists just of the
trivial group, and $\cC_n$ is defined inductively as the
class of groups which act on trees with all isotropy
subgroups belonging to $\cC_{n-1}$ \cite{Oy}.
\item The fundamental group of a closed manifold which is
  CAT(0)-cubical complex, i.e.~a piecewise Euclidean cubical complex
  which is 
non-positively curved in the sense of Alexandrov (cf.\
\eqref{def:cat})
\item Any word hyperbolic group in the sense of Gromov, in particular
  fundamental groups of strictly negatively curved manifolds (in the
  Riemannian sense or in the sense of Alexandrov).
\end{enumerate}
In the second case the claim that the assembly map is an
isomorphism follows from combining work of Niblo and Reeves
\cite{NR} and Higson and Kasparov \cite{HK}. Niblo and
Reeves \cite{NR}, compare \cite[p.\ 158]{Ju} show that any
group acting properly on a CAT(0)-cubical complex has the Haagerup
property  (in other words, is
a-T-menable).
Higson and Kasparov
prove that for such groups the Baum-Connes Conjecture holds.

The third case is proved by Mineyev and Yu \cite{MY} (reducing to some
deep results of Lafforgue \cite{L1,L2}). Note that the
fundamental group of a negatively curved manifold is known to be word
hyperbolic.

We note
that the papers cited above deal with the complex version of
the Baum-Connes Conjecture,  whereas we are interested in the
real version. However, it is a ``folk theorem'' that the 
complex isomorphism implies the real isomorphism, compare
e.g.~\cite{Ka}.
\end{num}

Unfortunately, making these kinds of assumptions about the
group $\pi$ does not save Conjectures \ref{num:GLconj} or
\ref{num:GLRconj}:

\begin{num}{\bf Addendum to Theorem \ref{thm:mainthm}.}
\label{num:addendum}
The group $\pi$ in Theorem \ref{thm:mainthm} may be chosen to
be in the class $\cC$ and to be the fundamental group of a
closed manifold which is a CAT(0)-cubical complex. Alternatively, we
may
choose
$\pi$ to be the fundamental group of a closed manifold which
is negatively curved in the Alexandrov sense.
\end{num}

With this `negative' result showing that Conjectures
\ref{num:GLconj} and \ref{num:GLRconj} are false even for
very ``reasonable" groups $\pi$, the challenge is to find a
necessary and sufficient condition for $M$ to
admit a \pscm. To the authors knowledge there is
that this point in time not even a {\it candidate} for this
condition.

\section{Outline of the proof}

As a first step we rephrase Conjectures \ref{num:GLconj} and
\ref{num:GLRconj}. The key for this is the following result
of Gromov-Lawson \cite{GL1} (see also \cite{RS}):
\begin{theorem}\label{thm:bordismthm}
 Let $M$ be a spin manifold
of dimension $n\ge 5$, and let
$u\colon M\to B\pi$ be the classifying map of its universal
covering. Then $M$ admits a \pscm\ if and only if $[M,u]\in
\Omega_n^{spin,+}(B\pi)$.
\end{theorem}
Here for any space $X$, the group 
$\Omega_n^{spin,+}(X)$ is by definition the following subgroup
of $\Omega_n^{spin}(X)$:
\[
\Omega_n^{spin,+}(X)=\left\{
[N,f]\Biggr|
\begin{gathered}
\text{$N$ is an $n$\nb-dimensional spin manifold}\\
\text{with \pscm,}\\
f\colon N\to X
\end{gathered}
\right\}
\] 
It
should be emphasized that while the `only if' portion of the
above theorem is of course tautological, the `if' part is
{\it not}; it can be rephrased by saying that if
$M$ is bordant (over $B\pi$) to some spin manifold $N$ with a \pscm,
then
$M$ itself admits a \pscm. 

We claim that for a finitely presented group $\pi$ the
Gromov-Lawson Conjecture \ref{num:GLconj} (resp.\
Gromov-Lawson-Rosenberg Conjecture \ref{num:GLRconj}) is
equivalent to the first (resp.\ second) of the following
equalities for $n\ge 5$
\begin{align}
\label{eq:GLconj} \Omega_n^{spin,+}(B\pi)&=\ker p\circ D\\
\label{eq:GLRconj} \Omega_n^{spin,+}(B\pi)&=\ker\alpha
\end{align}
By the Bordism Theorem \ref{thm:bordismthm} it is clear that
these equations imply Conjectures \ref{num:GLconj} resp.\
\ref{num:GLRconj}. Conversely, if $\pi$ is a finitely
presented group, then a surgery argument shows that {\it
every} bordism class in $\Omega_n^{spin}(B\pi)$, $n\ge 5$ can
be represented by a pair $(M,u)$, where $M$ is a manifold with
fundamental group $\pi$, and $u\colon M\to B\pi$ is the
classifying map of the universal covering of $M$. This shows
that the Gromov-Lawson Conjecture \ref{num:GLconj} implies
equation \eqref{eq:GLconj} and that the
Gromov-Lawson-Rosenberg Conjecture \ref{num:GLRconj} implies
equation \eqref{eq:GLRconj}.

Using this translation, our main result Theorem
\ref{thm:mainthm} is a consequence of the following slightly
more precise result.

\begin{theorem}\label{thm:maintechnical}
For $5\le n\le 8$, there  are finitely presented torsion free
groups $\pi$ such that 
$\Omega_n^{spin,+}(B\pi)\not\subseteq \ker p\circ D$.
\end{theorem}
To prove this statement, we need to produce a bordism class
$x\in\Omega_n^{spin}(B\pi)$ in the kernel of $p\circ D$ and
show that it cannot be represented by a manifold of \psc. At
present, three methods are known to show that a manifold
$M$ does not admit a
\pscm: the index-theory of the Dirac operator on
$M$, the Seiberg-Witten invariants, and the stable
minimal hypersurface method pioneered by Schoen and Yau
\cite{SY}. For the case at hand, the first two methods are
useless: the index of any manifold $M$ representing $x$
vanishes due to our
assumption
$p\circ D(x)=0$ and the Seiberg-Witten invariants of $M$ are
only defined for
$4$\nb-dimensional manifolds. 
As explained in  \cite[Proof of Cor.\ 1.5]{Sch}, a
corollary of the stable minimal hypersurface method is the
following result.

\begin{theorem}[\cite{Sch}, Cor.\ 1.5]
\label{thm:minhypersurface} 
Let $X$ be a space, and let
$H^+_n(X;\Z)$ be the subgroup of $H_n(X;\Z)$ consisting of
those elements which are of the form $f_*[N]$, where $N$ is
an oriented closed manifold of dimension $n$ which admits a \pscm, and
$f$ is a map $f\colon N\to X$. Let 
\[
\alpha\cap\colon H_n(X;\Z)\ra H_{n-1}(X;\Z)
\]
be the homomorphism given by the cap product with a class
$\alpha\in H^1(X;\Z)$. Then for $3\le n\le 8$ the homomorphism
$\alpha\cap$ maps $H_n^+(X;\Z)$ to $H_{n-1}^+(X;\Z)$.
\end{theorem} 

We remark that a better regularity result for hypersurfaces
with minimal volume representing a given homology class
proved by Smale \cite{Sm} makes it possible to include the 
case $n=8$ in the above result (see
\cite[\S  4]{JS} for a detailed  explanation).

\begin{cor}\label{cor:minhypersurface} Let $X$ be the
classifying space of a discrete group $\pi$. Assume that 
$x\in
\Omega_n^{spin}(X)$ is a bordism class for
$5\le n\le 8$, satisfying the condition 
\begin{equation}\label{eq:bordismcondition}
 \alpha_1\cap\dots\cap\alpha_{n-2}\cap
H(x)\ne 0\in H_2(X;\Z)
\end{equation}
for some cohomology classes 
$\alpha_1,\dots,\alpha_{n-2}\in H^1(X;\Z)$  
(here $H\colon \Omega_n^{spin}(X)\to H_n(X;\Z)$ is the
Hurewicz map given by sending a bordism class $[M,f]$ to
$f_*[M]$, where $[M]\in H_n(M;\Z)$ is the fundamental class of
$M$). Then
$x$ is not in the subgroup
$\Omega_n^{spin,+}(X)$.
\end{cor}

\begin{proof} Assume $x\in \Omega_n^{spin,+}(X)$. Then
$H(x)\in H_n^+(X;\Z)$ (by definition of $H_n^+(X;\Z)$);
applying Theorem \ref{thm:minhypersurface} first to the
cohomology class $\alpha_1$, then $\alpha_2$, e.t.c., we
conclude $\alpha_1\cap\dots\cap\alpha_{n-2}\cap
H(x)\in H_2^+(X;\Z)$. This is the desired contradiction,
since the element $\alpha_1\cap\dots\cap\alpha_{n-2}\cap
H(x)$ is assumed non-zero, while the group $H_2^+(X;\Z)$ is
trivial: by the Gauss-Bonnet Theorem, the only closed
oriented $2$\nb-manifold $N$ with \psc\ are disjoint unions
of $2$\nb-spheres; however any map $f$ from such a union to
the classifying space of a discrete group is homotopic to
the constant map, which implies $f_*[N]=0$.
\end{proof}

We note that Corollary \ref{cor:minhypersurface} implies 
Theorem
\ref{thm:maintechnical}, provided we can find a bordism class
$x$ in the kernel of $p\circ D$ satisfying condition
\ref{eq:bordismcondition}. Whether there is such a bordism
class is in general a pretty hard question; fortunately the
following homological condition is much easier to check, and,
as Theorem
\ref{thm:ssthm} below shows, it implies the
existence of a spin bordism class $x$ in the kernel of
$p\circ D$ with property 
 \eqref{eq:bordismcondition}.

\begin{num}{\bf Homological
condition.}\label{num:homcondition}
 There are  (co)homology classes
$\alpha_1,\dots,\alpha_{n-2}\in
H^1(X;\Z)$ and 
$z\in H_{n+5}(X;\Z)$  such that 
\begin{equation*}
\alpha_1\cap\dots\cap\alpha_{n-2}\cap\delta P^1\rho(z)\ne
0\in H_2(X;\Z)
\end{equation*}
Here 
\begin{itemize}
\item $\rho\colon H^*(X;\Z)\ra
H^*(X;\Z/3)$ is mod
$3$ reduction, 
\item $P^1\colon H_*(X;\Z/3)\lra H_{*-4}(X;\Z/3)$
is the homology operation dual to the degree
$4$ element $P^1$ of the mod
$3$ Steenrod algebra (see \cite[Chapter VI, section
1]{SE}), and
\item $\delta\colon H_*(X;\Z/3)\to H_{*-1}(X;\Z)$ is the
Bockstein homomorphism associated to the short exact
coefficient sequence $\Z\overset{\times 3}\ra
\Z\ra\Z/3$ (i.e., the boundary homomorphism of the
corresponding long exact homology sequence).
\end{itemize}
\end{num}

In section \ref{sec:ss} we will use the Atiyah-Hirzebruch
spectral sequences converging to $\Omega_*^{Spin}(X)$
resp.\ $KO_*(X)$ to prove the following result.

\begin{theorem}\label{thm:ssthm}
Let $X$ be a space which satisfies the
homological condition
\ref{num:homcondition} for $5\le n\le 8$.
Then there is a bordism
class $x\in
\Omega_n^{spin}(X)$ in the kernel of $p\circ D$ satisfying
condition \ref{eq:bordismcondition} above.
\end{theorem}

Putting these results together, we obtain the
following corollary.

\begin{cor}\label{cor:maincor} Let $X$ be a space
satisfying the cohomological condition \ref{num:homcondition}
for $5\le n\le 8$. If $X$ is the classifying space of
some discrete group $\pi$, then 
$\Omega_n^{spin,+}(B\pi)\not\subseteq \ker p\circ D$.
\end{cor}

 For a given space $X$ it is not hard
to check the homological condition \ref{num:homcondition},
provided we know enough about the (co)homology of $X$.
For example, in section \ref{example} we will prove the
following result by a straightforward calculation.

\begin{prop}\label{prop:example} For $n\ge 2$, let
$\Gamma_n$ be the cartesian product of
$2$ copies of $\Z/3$ and $n-2$
copies of $\Z$. Then the classifying space $B\Gamma_n$
satisfies the homological condition \ref{num:homcondition}.
\end{prop}

In particular,  
Corollary \ref{cor:maincor} then shows that the
Gromov-Lawson Conjecture
\ref{num:GLconj} does not hold for the group $\Gamma_n$ for
$5\le n\le 8$. This is very similar to the result of
one of the coauthors \cite{Sch} who
constructed a $5$\nb-dimensional spin manifold with
fundamental group $\pi=\Z/3\times \Z^4$ whose index
obstruction $\alpha([M,u])\in KO_5(C^*\pi)$ is
trivial, but which does not admit a metric of
\psc. However, the example in \cite{Sch} does {\it not} provide
a counter example to the Gromov-Lawson Conjecture
\ref{num:GLconj}, since it can be shown that $p\circ
D([M,u])\in KO_5(B\pi)$ is {\it non-trivial}. 

 As
explained in the  introduction, it is more interesting
to find  {\it torsion free} groups
$\pi$ for which the conjecture goes wrong. It seems
conceivable that experts might know explicit examples of
torsion free groups and enough about their
(co)homology to conclude that the homological
condition \ref{num:homcondition} is satisfied (or the
cohomological condition \ref{eq:cohomcondition}, which
is stronger, but easier to check). 

Lacking this
expertise, we argue more indirectly to prove Theorem
\ref{thm:maintechnical}. We use a construction of Baumslag,
Dyer and Heller \cite{BDH}, who associate a discrete group
$\pi$ to any connected CW complex $X$ and show that there is
a map $B\pi\to X$ which is an isomorphism in homology.
Moreover, if  $X$ is a finite CW complex, then $B\pi$ has the
homotopy type of a finite CW complex. In particular,  the
group $\pi$ is finitely presented. We would like to mention
that originally Kan and Thursten described a similar
construction \cite{KT}. However, their groups are usually
{\it not} finitely presented and hence not suitable for
our purposes. 

To prove Theorem \ref{thm:maintechnical} we let $X$ be any
finite CW complex satisfying the condition
\ref{num:homcondition}, e.g., the $n+5$\nb-skeleton of
$B\Gamma_n$ and let $\pi$ be the discrete group obtained from
$X$ via the Baumslag-Dyer-Heller procedure. Then $\pi$ is a
discrete group which is finitely presented and torsion free
since its classifying space $B\pi$ is homotopy equivalent to
a finite CW complex. Moreover, $B\pi$ satisfies the
homological condition \ref{num:homcondition} and hence
Corollary \ref{cor:maincor} implies Theorem
\ref{thm:maintechnical}.

Applying more sophisticated `asphericalization procedures'
due to Davis-Januskiewicz \cite{DJ} respectively Charney-Davis
\cite{CD}, we can produce groups $\pi$ satisfying the
geometric conditions mentioned in Addendum
\ref{num:addendum}. This is explained in section
\ref{sec:asphericalization}.


\section{A spectral sequence argument}\label{sec:ss}

The goal of this section is the proof of Theorem
\ref{thm:ssthm}. We consider the Atiyah-Hirzebruch
spectral sequence (AHSS for short)
\begin{equation}\label{AHSS}
E_{p,q}^2(X)=H_p(X;\Omega_q^{spin})\Longrightarrow
\Omega_{p+q}^{spin}(X).
\end{equation}
We recall that the Hurewicz map $H\colon
\Omega_n^{spin}(X)\lra H_n(X;\Z)$ is the {\it edge
homomorphism} of this spectral sequence; i.e., $H$ is
equal to the composition
\[
\Omega_n^{spin}(X)\onto E^\infty_{n,0}(X)\into
E^2_{n,0}(X)=H_n(X;\Z).
\]
This shows that in order to produce a bordism class
$x\in\Omega_n^{spin}(X)$ satisfying condition
\eqref{eq:bordismcondition} it suffices to produce a homology class
$y\in H_n(X;\Z)$ such that
\begin{gather}\label{property3}
y\in H_n(X;\Z)=E^2_{n,0}\quad\text{is an infinite
cycle for the AHSS \ref{AHSS}}\\
\label{property4}\alpha_1\cap\dots\cap\alpha_{n-2}\cap
y\ne 0
\end{gather}

To guarantee that the element $x$ produced this way is in the
kernel of the natural
transformation $p\circ D\colon \Omega_n^{Spin}(X)\ra
KO_n(X)$, we note that $p\circ D$  induces a map of spectral
sequences
\[
E^r_{p,q}(X;\Omega^{Spin})\lra E^r_{p,q}(X;KO).
\] 
To simplify the analysis of these spectral
sequences, from now on we localize all homology theories
at the prime $3$, which has the effect of replacing all
homology groups as well as the groups in the AHSS
converging to them by the corresponding localized
groups; i.e., their tensor product with
$\Z_{(3)}=\{\frac ab\mid \text{$b$ is prime to $3$}\}$. Note that we
continue to write e.g.~$H_*(X;\Z)$, but mean the localized group.

In particular, the coefficient ring $KO_*$ is now the
ring of polynomials $\Z_{(3)}[b]$ with generator $b\in
KO_4$. This implies that in the AHSS converging to
$KO_*(X)$ only the rows $E^r_{p,q}$ for $q\equiv
0\mod{4}$ are possibly non-trivial. In particular, the
first differential that can be non-trivial is $d_5$. To
finish the proof of Theorem \ref{thm:ssthm} we will need
the following result, which identifies
$d_5$ with a homology operation. This is certainly
well known among experts; since we  failed to find an
explicit reference in the literature, for completeness
we include a proof of this result below.

\begin{lem}\label{differential} In the AHSS
converging to
$KO_*(X)$ localized at $3$, the differential 
\[
d_5\colon E^5_{p,q}=E^2_{p,q}=H_{p}(X;\Z)
\lra E^5_{p-5,q+4}=E^2_{p-5,q+4}=H_{p-5}(X;\Z)
\]
for $q\equiv 0\mod{4}$ is the composition
\begin{equation}\label{homologyoperation}
H_p(X;\Z)\overset{\rho}\lra
H_p(X;\Z/3)\overset{P^1}\lra
H_{p-4}(X;\Z/3)\overset{\delta}\lra\
H_{p-5}(X;\Z).
\end{equation}
For the notation, compare \ref{num:homcondition}.
\end{lem}

\begin{proof}[Proof of Theorem \ref{thm:ssthm}] Assume that
the space $X$ satisfies the homological condition
\ref{num:homcondition}; i.e., there are (co)homology classes
$\alpha_1\dots,\alpha_{n-2}\in H^1(X;\Z)$ and $z\in
H_{n+5}(X;\Z)$ with 
\[
\alpha_1\cap\dots\cap\alpha_{n_2}\cap y\ne 0
\qquad\text{for}\qquad
y\=\delta\circ P^1\circ\rho(z).
\]
 We observe that $y$ has
property \ref{property4} by our assumption on $z$;
moreover, it also has property \ref{property3} (i.e.,
it survives to the $E^\infty$\nb-term of the AHSS
converging to $\Omega_*^{spin}(X)$) by the following
argument. The map of spectral sequences
\[
E^r_{p,q}(X;\Omega^{Spin})
\lra
E^r_{p,q}(X;KO)
\]
is an isomorphism on the rows $q=0,4$ for $r=2$ and
hence for $r=3,4,5$ (all the other rows in the range
$-4<q<8$ are trivial since we work localized at the
prime $3$). Lemma
\ref{differential} shows that the differential 
\[
d_5\colon E^5_{n+5,-4}(X;KO)\lra E^5_{n,0}(X;KO)
\]
sends $z\in
H_{n+5}(X;\Z)=E^2_{n+5,-4}(X;KO)=E^5_{n+5,-4}(X;KO)$
to  $y\in E^5_{n,0}(X;KO)$. In particular, $y$ is in the
kernel of $d_5$ in the spectral sequence converging to
$KO_*(X)$, and hence also in the kernel of $d_5$ in the spectral
sequence converging to $\Omega_*^{Spin}(X)$. 
The next possibly non-trivial differential 
\[
d_9\colon E^9_{n,0}(X;\Omega^{Spin})
\ra E^9_{n-9,8}(X;\Omega^{Spin})
\]
 is trivial due to our
assumption $n\le 8$. 

This shows that there is a bordism class
$x\in\Omega_n^{Spin}(X)$ with $H(x)=y$. Next we want to
show that with a careful choice of $x$ we can also
arrange for $p\circ D(x)=0$. We note that the relation
$d_5(z)=y$ in the spectral sequence converging to
$KO_*(X)$ implies that $p\circ D(x)$ is zero in
\[
E^\infty_{n,0}(X;KO)=F_nKO_n(X)/F_{n-1}KO_n(X).
\]
 This
does {\it not} imply that $p\circ D(x)$ is zero, only
that
$p\circ D(x)$ lies in the filtration $n-1$ subgroup 
$F_{n-1}KO_n(X)\subset KO_n(X)$. 

We want to show that replacing $x$ by $x'=x-x''$ for a
suitable
$x''\in F_{n-1}\Omega_n^{Spin}(X)$ produces an element
with the desired properties $H(x')=H(x)$ and
$p\circ D(x')=0$. We note that the first condition is satisfied
because $H$ sends elements of
\[
F_{n-1}\Omega_n^{spin}(X)=\im\left(
\Omega_n^{spin}(X^{(n-1)})\lra
\Omega_n^{spin}(X)\right)
\]
to zero, since  $H_n(X^{(n-1)};\Z)=0$. To obtain
$p\circ D(x')=0$, we need to be able choose $x''$ such
that
$p\circ D(x'')=p\circ D(x)$; in other words, we need to
show that the map
\[
p\circ D\colon F_{n-1}\Omega_n^{spin}(X)\ra
F_{n-1}KO_n(X)
\]
 is surjective for $n\le 8$. The argument is the
following.
The map
$\Omega_q^{spin}\ra KO_q$ is an isomorphism for $0\le q
<8$ and surjective for
$q=8$. It follows that the map of spectral sequences 
\[
E^r_{p,q}(X;\Omega^{Spin}) \lra
E^r_{p,q}(X;KO)
\]
is a surjection for $r=2$, $0\le q\le 8$. Since all the
differentials of the domain spectral sequence in
the range $0< q\le 8$,
$n=p+q\le 8$ are trivial, the above map is also
surjective for $r=3,\dots,\infty$ in that range. The
groups $E^\infty_{p,q}$, $n=p+q\le 8$, $q>0$, are the
associated graded groups for the filtered groups
$F_{n-1}\Omega_n^{spin}(X)$ (resp.\ $F_{n-1}KO_n(X)$).
This shows that the map  
$F_{n-1}\Omega_n^{spin}(X)\ra F_{n-1}KO_n(X)$ is
surjective for $n\le 8$  and finishes the proof of 
Theorem \ref{thm:ssthm}.
\end{proof}

\begin{proof}[Proof of Lemma \ref{differential}]
Multiplication by the periodicity element $b\in KO_4$
produces a homotopy equivalence of spectra
$\Sigma^4KO\cong KO$, which in turn induces an
isomorphism of spectral sequences $E^r_{p,q}(X;KO)\cong
E^r_{p,q+4}(X;KO)$. Hence it suffices to prove the
corresponding statement for the differential $d_5\colon
E^5_{p,0}(X;KO)\ra E^5_{p-5,4}(X;KO)$.

Given an integer $k$, let $KO\langle
k\rangle\ra KO$ be the {\it $(k-1)$\nb-connected
cover} of $KO$. Up to homotopy equivalence,
$KO\langle k\rangle$ is characterized by the
properties that $\pi_n(KO\langle k\rangle)=0$
for $n<k$ and that the induced map 
\[
\pi_n(KO\langle k\rangle)\lra
\pi_n(KO)
\]
is an isomorphism for $n\ge k$. The spectrum
$KO\langle 0\rangle$ is also known as the
{\it connective real $K$\nb-theory spectrum\/} 
and is usually denoted $ko$. Given a second
integer $l\ge k$, let $KO\langle k,l\rangle$ be the
part of the Postnikov tower for $KO$, whose homotopy
groups $\pi_n(KO\langle k,l\rangle)$ are trivial for
$n<k$  or $n>l$ and are isomorphic to $\pi_n(KO)$ for
$k\le n\le l$ (this isomorphism is induced by a map
$KO\langle k\rangle\ra KO\langle k,l\rangle$).
In particular, $KO\langle k,k\rangle$ has only one
possibly non-trivial homotopy group and hence can be
identified with the Eilenberg-MacLane spectrum
$\Sigma^kH\pi_k(KO)$ (here $HA$ for an
abelian group $A$ is the Eilenberg-MacLane
spectrum characterized by $\pi_0(HA)=A$ and
$\pi_n(HA)=0$ for $n\ne 0$). 

We note that  the maps
$KO\langle 0\rangle\ra KO$ and $KO\langle 0\rangle\ra
KO\langle 0,4\rangle$ induce isomorphisms of AHSS's
in the range $0\le q\le 4$ for $r=2,3,4,5$ (since  we
work localized at $3$, all rows for $q\ne 0\mod{4}$
are trivial, and hence the first possibly non-trivial
differential is $d_5$). The AHSS
$E^r_{p,q}(X;KO\langle 0,4\rangle)$
has only two non-trivial rows and hence degenerates to
a long exact sequence
\begin{multline*}
\dots\lra  KO\langle 0,4\rangle_n(X)\lra
E^5_{n,0}=H_n(X;\Z)\overset{d_5}\lra\\
E^5_{n-5,4}=H_{n-5}(X;\Z)\lra 
KO\langle 0,4\rangle_{n-1}(X)\lra\dots .
\end{multline*}
This can be identified with the long exact homotopy
sequence induced by the Puppe sequence 
\begin{multline*}
KO\langle 4,4\rangle\wedge X\lra
KO\langle 0,4\rangle\wedge X\lra
KO\langle 0,0\rangle\wedge X=H\Z\wedge X\\
\overset{f\wedge 1}\lra
\Sigma KO\langle 4,4\rangle\wedge X=\Sigma^5 H\Z\wedge
X
\end{multline*}
The homotopy class of $f\colon H\Z\ra \Sigma^5H\Z$ can
be interpreted as a cohomology class in
$H^5(H\Z;\Z)\cong \Z/3$ (we've localized at $3$). The
generator of $H^5(H\Z;\Z)$ is given by applying
$\delta\circ P^1\circ\rho$ to the generator of
$H^0(H\Z;\Z)\cong\Z_{(3)}$. This shows that $d_5$ is a
{\it multiple} of the homology operation
\ref{homologyoperation}. 

To shows that it is a {\it non-trivial} multiple, it
suffices to show that the differential is non-trivial
for some space or spectrum $X$. We choose $X=H\Z/3$,
and consider the AHSS
\[
E^r_{p,q}(H\Z/3;ko)\Longrightarrow 
ko_{p+q}(H\Z/3)=\pi_{p+q}(ko\wedge
H\Z/3)=H_{p+q}(ko;\Z/3).
\]
It is well known that $H^*(ko;\Z/3)\cong
A/\left(AQ_1+A\beta\right)$. For this, see \cite[Prop.\ 2.3]{AP}; that
proposition applies to the spectrum $X=KO\langle 4\rangle$,
which localized at $3$ by Bott-periodicity may be identified with
$\Sigma^4ko$). Here
$A$ is the mod
$3$ Steenrod algebra, $\beta\in A_1$ is the
Bockstein ($A_n\subset A$ consists of the elements of
degree $n$), and
$Q_1\in A_5$ is the commutator of $P^1\in A_4$
and $\beta$. In particular, since $\beta P^1$
and $P^1\beta$ form a basis of $A_5$, the
cohomology group $H^5(ko;\Z/3)$ and hence the 
group $H_5(ko;\Z/3)=ko_5(H\Z/3)$ is trivial. It follows
that the non-trivial elements in 
\[
E^2_{5,0}(H\Z/3;ko)=H_5(H\Z/3;\Z)\cong\Z/3
\]
do not survive to the $E^\infty$\nb-term. For
dimensional reasons the only possibly non-trivial
differential is $d_5$. This finishes the proof of Lemma
\ref{differential}.
\end{proof}
\section{Construction of spaces satisfying the 
homological condition \ref{num:homcondition}}\label{example}

The goal of this section is the construction of spaces
satisfying the homological condition \ref{num:homcondition};
in particular, we will prove Proposition \ref{prop:example},
which claims that the classifying space
$B\Gamma_n$ satisfies this condition, where $\Gamma_n$ is the
product of two copies of $\Z/3$ and $n-2$ copies of $\Z$. For
this calculation it is convenient to pass to cohomology.  The
following lemma will give a {\it cohomological} condition which
implies the homological condition \ref{num:homcondition}. 

To
state the lemma, we first need some notation. Let 
$\rho\colon H^*(X;\Z)\to  H^*(X;\Z/3)$ be mod
$3$ reduction and let $\beta\colon H^*(X;\Z/3)\ra
H^{*+1}(X;\Z/3)$ the mod $3$ Bockstein homomorphism, the
boundary map of the long exact cohomology sequence
induced by the short exact coefficient sequence
$\Z/3\overset{\times 3}\lra\Z/3^2\lra\Z/3$. We recall
that $\beta$ is the composition $\rho\delta$, where $\delta\colon
H^*(X;\Z/3)\to H^{*+1}(X;\Z)$ is the integral Bockstein. Let
$Q_1$ be the degree $5$ element of the mod $3$ Steenrod
algebra which is the commutator $Q_1=[P^1,\beta]$.

\begin{lem} Let $X$ be a space and 
assume there are cohomology classes
$\alpha_1,\dots,\alpha_{n-2}\in H^1(X;\Z)$
and $\zeta\in H^2(X;\Z/3)$ such that
\begin{equation}\label{eq:cohomcondition}
\rho(\alpha_1)\cup\dots\cup\rho(\alpha_{n-2})\cup
\beta Q_1\zeta\ne 0\in H^{n+6}(X;\Z/3).
\end{equation}
Then $X$ satisfies the homological condition
\ref{num:homcondition}.
\end{lem}

\begin{proof} The assumption of the lemma implies that
there is a homology class $y\in H_{n+6}(X;\Z/3)$ such
that the Kronecker product
\[
\langle
\rho(\alpha_1)\cup\dots\cup\rho(\alpha_{n-2})\cup
\beta Q_1\zeta,
y\rangle
\]
is non-zero. We calculate, using that $Q_1$ and $\beta$ are graded
derivations and that the cohomological $\beta$ and $Q_1$ are dual to the
homological versions,
\begin{align*}
 0& \ne \langle
\rho(\alpha_1)\cup\dots\cup\rho(\alpha_{n-2})\cup
\beta Q_1\zeta,
y\rangle\\
&=\langle
\rho(\alpha_1)\cup\dots\cup\rho(\alpha_{n-2})\cup
Q_1\zeta,
\beta y\rangle\qquad\text{(since $\beta\rho=0$)}\\
&=\langle
\rho(\alpha_1)\cup\dots\cup\rho(\alpha_{n-2})\cup
\zeta,
Q_1(\beta y)\rangle\qquad\text{(since $Q_1(H^1(X))=0$)}\\
&=\langle
\rho(\alpha_1)\cup\dots\cup\rho(\alpha_{n-2})\cup
\zeta,
\beta P^1 \rho(\delta y)\rangle\qquad\text{(since $\beta\beta=0$)}\\
&=\langle\zeta,
\rho\left(\alpha_1\cap\dots\cap\alpha_{n-2}\cap
\delta P^1\rho(\delta y)\right)\rangle\qquad\text{(since $\beta=\rho\delta$)}
\end{align*}
This shows that the
non-triviality of $\rho(\alpha_1)\cup\dots\cup\rho(\alpha_{n-2})\cup
\beta Q_1\zeta$ implies that 
$\alpha_1\cap\dots\cap\alpha_{n-2}\cap
\delta P^1\rho(z)$ for $z=\delta y$ is non-trivial. 
\end{proof}

Recall the statement of Proposition \ref{prop:example}:
\begin{proprecall} For $n\ge 2$, let
$\Gamma_n$ be the Cartesian product of
$2$ copies of $\Z/3$ and $n-2$
copies of $\Z$. Then the classifying space $B\Gamma_n$
satisfies the homological condition \ref{num:homcondition}.
\end{proprecall}

\begin{proof}[Proof of Proposition \ref{prop:example}] We recall
that the cohomology ring of 
\[
B\Gamma_n=\underbrace{S^1\times\dots\times S^1}_{n-2}\times
B\Z/3\times B\Z/3
\]
is given by 
\[
H^*(B\Gamma_n;\Z/3)=
\Lambda\left(\rho(\alpha_1),\dots,\rho(\alpha_{n_2}),x_1,x_2\right)
\otimes \Z/3[\beta x_1,\beta x_2],
\]
where $\alpha_i\in H^1(B\Gamma_n;\Z)$ is the pull back of the
generator of $H^1(S^1;\Z)$ via the projection to the $i$-th copy
of $S^1$, and $x_1$ (resp.\ $x_2$) is  the pull
back of  the generator of $H^1(B\Z/3;\Z/3)$ via the
projection to the first (resp.\ second) $B\Z/3$\nb-factor
of $B\Gamma_n$.  

We calculate for $j=1,2$
\[
Q_1(x_j)=[P^1,\beta]x_j=P^1\beta x_j-\beta P^1
x_j=(\beta x_j)^3.
\]
Since $\beta$ and $Q_1$ are graded derivations, it
follows that 
\[
\beta Q_1(x_1x_2)=\beta\left((\beta x_1)^3x_2
-x_1(\beta x_2)^3\right)
=(\beta x_1)^3\beta x_2-\beta x_1(\beta
x_2)^3\ne 0.
\]
It follows that $\zeta=x_1x_2$ satisfies the cohomological
condition \eqref{eq:cohomcondition}.
\end{proof}
\section{Asphericalization
procedures}\label{sec:asphericalization}

In this section we prove Addendum \ref{num:addendum} to
our main theorem claiming that the groups $\pi$ for which we
can construct counter examples to the Gromov-Lawson Conjecture
\ref{num:GLconj} may be chosen to have classifying spaces
which are manifolds which are non-positively (resp.\
negatively) curved in the Alexandrov sense. This means the
following.

\begin{defn}\label{def:cat}
  A \emph{length space} is a metric space where any two points can be
  joined by a geodesic. A length space $X$ is called a (locally)
  \emph{CAT($r$)-space} for $r\in\R$, if given a triangle $T$
in a
  sufficiently small open set of $X$ and a vertex $x$ of $T$, the
  distance from $x$ to the opposite edge is not more than the
  corresponding distance in a comparison triangle in a simply
  connected manifold with constant curvature equal to $r$. A
  \emph{comparison triangle} is a triangle with the same side lengths
  as the given one.

  We say that $X$ is \emph{negatively curved in the
Alexandrov sense} (resp.\ \emph{non-positively
  curved}), if $X$ is a CAT($r$)-space
  with $r<0$ (resp.\ $r=0$).
\end{defn}
This is a generalization of the classical notion: every complete
Riemannian manifold of negative curvature is negatively curved in the
Alexandrov sense, and correspondingly for non-positive curvature.

To produce discrete groups $\pi$ whose classifying space is a
manifold which is non-positively curved (resp.\ negatively curved)
in the Alexandrov sense, we use asphericalization procedures due
to Davis-Januskiewicz \cite{DJ} and Charney-Davis \cite{CD}. 
An {\it asphericalization procedure} assigns to every space $X$ in
a certain class of spaces an aspherical space $BLX$
together with a map $BLX\to X$ with certain homological
properties. Here the notation $BLX$ is chosen to reflect the fact
that $BLX$ serves as a classifying space for its fundamental
group, for which the notation $LX$ is used.

Typically, the spaces $X$ considered are simplicial (or
cubical) cell complexes and $BLX$ is constructed by replacing each
$n$\nb-cell of $X$ by some `model space' and by gluing together
these model spaces to form $BLX$ according to the same
combinatorial patters as the cells of $X$ are glued to form $X$.

If this is done carefully, as in the procedures described below, one
can construct metrics on the result which satisfy appropriate curvature
conditions. However, even if the resulting space is a smooth manifold,
it is not clear at all whether the metric can be chosen to be a smooth
Riemannian metric.

\subsubsection*{Baumslag-Dyer-Heller asphericalization}
The construction of Baumslag, Dyer and Heller uses as basic
building block an acyclic group, i.e.~a non-trivial group with trivial
integral homology. For suitable choices of this
building block, their construction gives explicit descriptions of 
$LX$ as elements of $\cC$. This description uses the
combinatorics of the simplicial complex $X$ (compare also the
description of a very similar asphericalization procedure in
\cite{Bl}). The map $BLX\to X$ they produce induces an isomorphism 
in homology.

\subsubsection*{Davis-Januskiewicz asphericalization}

The goal of constructing $BLX$ which is non-positively curved can be
achieved using an asphericalization procedure of Davis and
Januskiewicz \cite{DJ}. If $X$ is a closed $n$\nb-manifold,
their construction produces a new $n$\nb-manifold $BLX$ which
is non-positively curved in the Alexandrov sense, and whose
fundamental group $LX$ belongs to the class $\cC$. The price to be
paid is that  the  map
$BLX\to X$ unlike in the case of the Baumslag-Dyer-Heller
procedure in general does {\it not} induces an isomorphism in
homology (necessarily so: e.g.\ for
$X=S^2$ there is no non-positively curved $2$\nb-manifold with the
same homology as $S^2$). However, the map $H_*(BLX;\Z)\to
H_*(X;\Z)$ is still {\it surjective}.

This is good enough to prove part (1) of Addendum
\ref{num:addendum} as follows. Let $Y$ be the manifold
with boundary obtained as a `thickening' of the
$n+5$\nb-skeleton of
$B\Gamma_n$ (cf.\ Prop.\ \ref{prop:example}) in $\R^N$ with $N$
large. Let $X$ be the boundary of $Y$
We note that the homology of $X$ is isomorphic to the
homology of $B\Gamma_n$ in degrees $\le n+5$ for $N$ sufficiently
large. In particular, since $B\Gamma_n$ satisfies the homology
condition \ref{num:homcondition}, so does the manifold
$X$.  Then the surjectivity of the map $H_*(BLX;\Z)\to H_*(X;\Z)$
implies that also $BLX$ satisfies the condition, and hence
Corollary \ref{cor:maincor} implies that the Gromov-Lawson
Conjecture does not hold for $\pi=LX$.

\subsubsection*{Charney-Davis asphericalization}
This procedure is essentially a strengthening of  the
Davis-Januskiewicz procedure: from a closed $n$\nb-manifold $X$
it produces an $n$\nb-manifold $BLX$ which is  {\it negatively
curved} in the Alexandrov sense and a map $BLX\to X$ which induces
a surjection on homology. However, $LX$ might not belong to the
class $\cC$. 

With the same argument as above, we can then produce a group $\pi$
for which the  Gromov-Lawson Conjecture doesn't hold, which is the
fundamental group of a negatively curved manifold. This proves the
second part of Addendum \ref{num:addendum}.


\end{document}